
\input amstex
\expandafter\ifx\csname mathdefs.tex\endcsname\relax
  \expandafter\gdef\csname mathdefs.tex\endcsname{}
\else \message{Hey!  Apparently you were trying to
  \string twice.   This does not make sense.} 
\errmessage{Please edit your file (probably \jobname.tex) and remove
any duplicate ``\string\input'' lines} \fi




\catcode`\X=12\catcode`\@=11

\def\n@wcount{\alloc@0\count\countdef\insc@unt}
\def\n@wwrite{\alloc@7\write\chardef\sixt@@n}
\def\n@wread{\alloc@6\read\chardef\sixt@@n}
\def\r@s@t{\relax}\def\v@idline{\par}\def\@mputate#1/{#1}
\def\l@c@l#1X{\firstpart.#1}\def\gl@b@l#1X{#1}\def\t@d@l#1X{{}}

\def\crossrefs#1{\ifx\all#1\let\tr@ce=\all\else\def\tr@ce{#1,}\fi
   \n@wwrite\cit@tionsout\openout\cit@tionsout=\jobname.cit 
   \write\cit@tionsout{\tr@ce}\expandafter\setfl@gs\tr@ce,}
\def\setfl@gs#1,{\def\@{#1}\ifx\@\empty\let\next=\relax
   \else\let\next=\setfl@gs\expandafter\xdef
   \csname#1tr@cetrue\endcsname{}\fi\next}
\def\m@ketag#1#2{\expandafter\n@wcount\csname#2tagno\endcsname
     \csname#2tagno\endcsname=0\let\tail=\all\xdef\all{\tail#2,}
   \ifx#1\l@c@l\let\tail=\r@s@t\xdef\r@s@t{\csname#2tagno\endcsname=0\tail}\fi
   \expandafter\gdef\csname#2cite\endcsname##1{\expandafter
     \ifx\csname#2tag##1\endcsname\relax?\else\csname#2tag##1\endcsname\fi
     \expandafter\ifx\csname#2tr@cetrue\endcsname\relax\else
     \write\cit@tionsout{#2tag ##1 cited on page \folio.}\fi}
   \expandafter\gdef\csname#2page\endcsname##1{\expandafter
     \ifx\csname#2page##1\endcsname\relax?\else\csname#2page##1\endcsname\fi
     \expandafter\ifx\csname#2tr@cetrue\endcsname\relax\else
     \write\cit@tionsout{#2tag ##1 cited on page \folio.}\fi}
   \expandafter\gdef\csname#2tag\endcsname##1{\expandafter
      \ifx\csname#2check##1\endcsname\relax
      \expandafter\xdef\csname#2check##1\endcsname{}%
      \else\immediate\write16{Warning: #2tag ##1 used more than once.}\fi
      \multit@g{#1}{#2}##1/X%
      \write\t@gsout{#2tag ##1 assigned number \csname#2tag##1\endcsname\space
      on page \number\count0.}%
   \csname#2tag##1\endcsname}}
\def\multit@g#1#2#3/#4X{\def\t@mp{#4}\ifx\t@mp\empty%
      \global\advance\csname#2tagno\endcsname by 1 
      \expandafter\xdef\csname#2tag#3\endcsname
      {#1\number\csname#2tagno\endcsnameX}%
   \else\expandafter\ifx\csname#2last#3\endcsname\relax
      \expandafter\n@wcount\csname#2last#3\endcsname
      \global\advance\csname#2tagno\endcsname by 1 
      \expandafter\xdef\csname#2tag#3\endcsname
      {#1\number\csname#2tagno\endcsnameX}
      \write\t@gsout{#2tag #3 assigned number \csname#2tag#3\endcsname\space
      on page \number\count0.}\fi
   \global\advance\csname#2last#3\endcsname by 1
   \def\t@mp{\expandafter\xdef\csname#2tag#3/}%
   \expandafter\t@mp\@mputate#4\endcsname
   {\csname#2tag#3\endcsname\lastpart{\csname#2last#3\endcsname}}\fi}
\def\t@gs#1{\def\all{}\m@ketag#1e\m@ketag#1s\m@ketag\t@d@l p
   \m@ketag\gl@b@l r \n@wread\t@gsin
   \openin\t@gsin=\jobname.tgs \re@der \closein\t@gsin
   \n@wwrite\t@gsout\openout\t@gsout=\jobname.tgs }
\outer\def\localtags{\t@gs\l@c@l}
\outer\def\globaltags{\t@gs\gl@b@l}
\outer\def\newlocaltag#1{\m@ketag\l@c@l{#1}}
\outer\def\newglobaltag#1{\m@ketag\gl@b@l{#1}}

\newif\ifpr@ 
\def\m@kecs #1tag #2 assigned number #3 on page #4.%
   {\expandafter\gdef\csname#1tag#2\endcsname{#3}
   \expandafter\gdef\csname#1page#2\endcsname{#4}
   \ifpr@\expandafter\xdef\csname#1check#2\endcsname{}\fi}
\def\re@der{\ifeof\t@gsin\let\next=\relax\else
   \read\t@gsin to\t@gline\ifx\t@gline\v@idline\else
   \expandafter\m@kecs \t@gline\fi\let \next=\re@der\fi\next}
\def\pretags#1{\pr@true\pret@gs#1,,}
\def\pret@gs#1,{\def\@{#1}\ifx\@\empty\let\n@xtfile=\relax
   \else\let\n@xtfile=\pret@gs \openin\t@gsin=#1.tgs \message{#1} \re@der 
   \closein\t@gsin\fi \n@xtfile}

\newcount\sectno\sectno=0\newcount\subsectno\subsectno=0
\newif\ifultr@local \def\ultralocal{\ultr@localtrue}
\def\firstpart{\number\sectno}
\def\lastpart#1{\ifcase#1 \or a\or b\or c\or d\or e\or f\or g\or h\or 
   i\or k\or l\or m\or n\or o\or p\or q\or r\or s\or t\or u\or v\or w\or 
   x\or y\or z \fi}

\def\resetall{\global\advance\sectno by 1\subsectno=0
   \gdef\firstpart{\number\sectno}\r@s@t}
\def\resetsub{\global\advance\subsectno by 1
   \gdef\firstpart{\number\sectno.\number\subsectno}\r@s@t}
\def\newsection#1\par{\resetall\vskip0pt plus.3\vsize\penalty-250
   \vskip0pt plus-.3\vsize\bigskip\bigskip
   \message{#1}\leftline{\bf#1}\nobreak\bigskip}
\def\subsection#1\par{\ifultr@local\resetsub\fi
   \vskip0pt plus.2\vsize\penalty-250\vskip0pt plus-.2\vsize
   \bigskip\smallskip\message{#1}\leftline{\bf#1}\nobreak\medskip}

\def\t@gsoff#1,{\def\@{#1}\ifx\@\empty\let\next=\relax\else\let\next=\t@gsoff
   \def\@@{p}\ifx\@\@@\else
   \expandafter\gdef\csname#1cite\endcsname##1{\zeigen{##1}}
   \expandafter\gdef\csname#1page\endcsname##1{?}
   \expandafter\gdef\csname#1tag\endcsname##1{\zeigen{##1}}\fi\fi\next}
\def\verbatimtags{\ifx\all\relax\else\expandafter\t@gsoff\all,\fi}
\def\zeigen#1{\hbox{$\langle$}#1\hbox{$\rangle$}}

\def\(#1){\edef\dot@g{\ifmmode\ifinner(\hbox{\noexpand\etag{#1}})
   \else\noexpand\eqno(\hbox{\noexpand\etag{#1}})\fi
   \else(\noexpand\ecite{#1})\fi}\dot@g}

\newif\ifbr@ck
\def\eat#1{}
\def\[#1]{\br@cktrue[\br@cket#1'X]}
\def\br@cket#1'#2X{\def\temp{#2}\ifx\temp\empty\let\next\eat
   \else\let\next\br@cket\fi
   \ifbr@ck\br@ckfalse\br@ck@t#1,X\else\br@cktrue#1\fi\next#2X}
\def\br@ck@t#1,#2X{\def\temp{#2}\ifx\temp\empty\let\neext\eat
   \else\let\neext\br@ck@t\def\temp{,}\fi
   \def\teemp{#1}\ifx\teemp\empty\else\rcite{#1}\fi\temp\neext#2X}
\def\resetbr@cket{\gdef\[##1]{[\rtag{##1}]}}
\def\references{\resetbr@cket\newsection References\par}

\newtoks\symb@ls\newtoks\s@mb@ls\newtoks\p@gelist\n@wcount\ftn@mber
    \ftn@mber=1\newif\ifftn@mbers\ftn@mbersfalse\newif\ifbyp@ge\byp@gefalse
\def\defm@rk{\ifftn@mbers\n@mberm@rk\else\symb@lm@rk\fi}
\def\n@mberm@rk{\xdef\m@rk{{\the\ftn@mber}}%
    \global\advance\ftn@mber by 1 }
\def\rot@te#1{\let\temp=#1\global#1=\expandafter\r@t@te\the\temp,X}
\def\r@t@te#1,#2X{{#2#1}\xdef\m@rk{{#1}}}
\def\b@@st#1{{$^{#1}$}}\def\str@p#1{#1}
\def\symb@lm@rk{\ifbyp@ge\rot@te\p@gelist\ifnum\expandafter\str@p\m@rk=1 
    \s@mb@ls=\symb@ls\fi\write\f@nsout{\number\count0}\fi \rot@te\s@mb@ls}
\def\byp@ge{\byp@getrue\n@wwrite\f@nsin\openin\f@nsin=\jobname.fns 
    \n@wcount\currentp@ge\currentp@ge=0\p@gelist={0}
    \re@dfns\closein\f@nsin\rot@te\p@gelist
    \n@wread\f@nsout\openout\f@nsout=\jobname.fns }
\def\m@kelist#1X#2{{#1,#2}}
\def\re@dfns{\ifeof\f@nsin\let\next=\relax\else\read\f@nsin to \f@nline
    \ifx\f@nline\v@idline\else\let\t@mplist=\p@gelist
    \ifnum\currentp@ge=\f@nline
    \global\p@gelist=\expandafter\m@kelist\the\t@mplistX0
    \else\currentp@ge=\f@nline
    \global\p@gelist=\expandafter\m@kelist\the\t@mplistX1\fi\fi
    \let\next=\re@dfns\fi\next}
\def\symbols#1{\symb@ls={#1}\s@mb@ls=\symb@ls} 
\def\bigsymbol{\textstyle}
\symbols{\bigsymbol\ast,\dagger,\ddagger,\sharp,\flat,\natural,\star}
\def\ftnumbers{\ftn@mberstrue} \def\ftsymbols{\ftn@mbersfalse}
\def\paginal{\byp@ge} \def\resetftnumbers{\ftn@mber=1}
\def\ftnote#1{\defm@rk\expandafter\expandafter\expandafter\footnote
    \expandafter\b@@st\m@rk{#1}}

\long\def\jump#1\endjump{}
\def\ssum{\mathop{\lower .1em\hbox{$\textstyle\Sigma$}}\nolimits}

\def\qed{\nobreak\kern 1em \vrule height .5em width .5em depth 0em}
\def\newneq{\hbox{\rlap{\hbox to 1\wd9{\hss$=$\hss}}\raise .1em 
   \hbox to 1\wd9{\hss$\scriptscriptstyle/$\hss}}}
\def\subsetne{\setbox9 = \hbox{$\subset$}\mathrel{\hbox{\rlap
   {\lower .4em \newneq}\raise .13em \hbox{$\subset$}}}}
\def\supsetne{\setbox9 = \hbox{$\subset$}\mathrel{\hbox{\rlap
   {\lower .4em \newneq}\raise .13em \hbox{$\supset$}}}}

\def\vbar{\mathchoice{\vrule height6.3ptdepth-.5ptwidth.8pt\kern-.8pt}
   {\vrule height6.3ptdepth-.5ptwidth.8pt\kern-.8pt}
   {\vrule height4.1ptdepth-.35ptwidth.6pt\kern-.6pt}
   {\vrule height3.1ptdepth-.25ptwidth.5pt\kern-.5pt}}
\def\f@dge{\mathchoice{}{}{\mkern.5mu}{\mkern.8mu}}
\def\b@c#1#2{{\rm \mkern#2mu\vbar\mkern-#2mu#1}}
\def\b@b#1{{\rm I\mkern-3.5mu #1}}
\def\b@a#1#2{{\rm #1\mkern-#2mu\f@dge #1}}
\def\bb#1{{\count4=`#1 \advance\count4by-64 \ifcase\count4\or\b@a A{11.5}\or
   \b@b B\or\b@c C{5}\or\b@b D\or\b@b E\or\b@b F \or\b@c G{5}\or\b@b H\or
   \b@b I\or\b@c J{3}\or\b@b K\or\b@b L \or\b@b M\or\b@b N\or\b@c O{5} \or
   \b@b P\or\b@c Q{5}\or\b@b R\or\b@a S{8}\or\b@a T{10.5}\or\b@c U{5}\or
   \b@a V{12}\or\b@a W{16.5}\or\b@a X{11}\or\b@a Y{11.7}\or\b@a Z{7.5}\fi}}

\catcode`\X=11 \catcode`\@=12

\localtags
\sectno=-1   
\NoBlackBoxes
\documentstyle {amsppt}
\topmatter
\title\nofrills {In the random graph $G(n,p),p=n^{-a}$: \\
 if $\psi$ has
probability $O(n^{-\varepsilon})$ for every $\varepsilon > 0$ then
 it has probability $O(e^{-n^\varepsilon})$ for some $\varepsilon > 0$ \\
\bigskip
Sh551} \endtitle
\rightheadtext{Random Graph}
\bigskip
\author {Saharon Shelah \thanks{\null\newline
We thank Alice Leonhardt for the beautiful typing. \null\newline
Latest Revision 95/Dec/22 \null\newline
The author would like to thank the United States Israel Binational Science
Foundation for partially supporting this research. \newline
We thank Joel Spencer for telling us the problem} \endthanks} \endauthor
\affil {Institute of Mathematics \\
The Hebrew University
Jerusalem, Israel
\medskip
Rutgers University \\
Department of Mathematics \\
New Brunswick, NJ USA} \endaffil
\endtopmatter
\document  
\expandafter\ifx\csname bib4plain.tex\endcsname\relax
  \expandafter\gdef\csname bib4plain.tex\endcsname{}
\else \message{Hey!  Apparently you were trying to \string twice.   This does not make sense.}
\errmessage{Please edit your file (probably \jobname.tex) and remove
any duplicate ``\string\input'' lines} \fi

\def\renewcommand{\newcommand}	       
\edef\cite{\the\catcode`@}%
\catcode`@ = 11
\let\@oldatcatcode = \cite
\chardef\@letter = 11
\chardef\@other = 12
%
%
%
%
\def\@innerdef#1#2{\edef#1{\expandafter\noexpand\csname #2\endcsname}}%
%
%
\@innerdef\@innernewcount{newcount}%
\@innerdef\@innernewdimen{newdimen}%
\@innerdef\@innernewif{newif}%
\@innerdef\@innernewwrite{newwrite}%
%
%
%
\def\@gobble#1{}%
%
%
%
\ifx\inputlineno\@undefined
   \let\@linenumber = \empty 
\else
   \def\@linenumber{\the\inputlineno:\space}%
\fi
%
%
%
\def\@futurenonspacelet#1{\def\cs{#1}%
   \afterassignment\@stepone\let\@nexttoken=
}%
\begingroup 
\def\\{\global\let\@stoken= }%
\\ 
\endgroup
\def\@stepone{\expandafter\futurelet\cs\@steptwo}%
\def\@steptwo{\expandafter\ifx\cs\@stoken\let\@@next=\@stepthree
   \else\let\@@next=\@nexttoken\fi \@@next}%
\def\@stepthree{\afterassignment\@stepone\let\@@next= }%
%
%
%
\def\@getoptionalarg#1{%
   \let\@optionaltemp = #1%
   \let\@optionalnext = \relax
   \@futurenonspacelet\@optionalnext\@bracketcheck
}%
%
%
\def\@bracketcheck{%
   \ifx [\@optionalnext
      \expandafter\@@getoptionalarg
   \else
      \let\@optionalarg = \empty
      \expandafter\@optionaltemp
   \fi
}%
\def\@@getoptionalarg[#1]{%
   \def\@optionalarg{#1}%
   \@optionaltemp
}%
%
%
%
\def\@nnil{\@nil}%
\def\@fornoop#1\@@#2#3{}%
\def\@for#1:=#2\do#3{%
   \edef\@fortmp{#2}%
   \ifx\@fortmp\empty \else
      \expandafter\@forloop#2,\@nil,\@nil\@@#1{#3}%
   \fi
}%
\def\@forloop#1,#2,#3\@@#4#5{\def#4{#1}\ifx #4\@nnil \else
       #5\def#4{#2}\ifx #4\@nnil \else#5\@iforloop #3\@@#4{#5}\fi\fi
}%
\def\@iforloop#1,#2\@@#3#4{\def#3{#1}\ifx #3\@nnil
       \let\@nextwhile=\@fornoop \else
      #4\relax\let\@nextwhile=\@iforloop\fi\@nextwhile#2\@@#3{#4}%
}%
%
%
%
\@innernewif\if@fileexists
\def\@testfileexistence{\@getoptionalarg\@finishtestfileexistence}%
\def\@finishtestfileexistence#1{%
   \begingroup
      \def\extension{#1}%
      \immediate\openin0 =
         \ifx\@optionalarg\empty\jobname\else\@optionalarg\fi
         \ifx\extension\empty \else .#1\fi
         \space
      \ifeof 0
         \global\@fileexistsfalse
      \else
         \global\@fileexiststrue
      \fi
      \immediate\closein0
   \endgroup
}%
%
%
%
%
\def\bibliographystyle#1{%
   \@readauxfile
   \@writeaux{\string\bibstyle{#1}}%
}%
\let\bibstyle = \@gobble
%
%
\let\bblfilebasename = \jobname
\def\bibliography#1{%
   \@readauxfile
   \@writeaux{\string\bibdata{#1}}%
   \@testfileexistence[\bblfilebasename]{bbl}%
   \if@fileexists
      \nobreak
      \@readbblfile
   \fi
}%
\let\bibdata = \@gobble
%
%
\def\nocite#1{%
   \@readauxfile
   \@writeaux{\string\citation{#1}}%
}%
\@innernewif\if@notfirstcitation
%
%
\def\cite{\@getoptionalarg\@cite}%
%
%
\def\@cite#1{%
   \let\@citenotetext = \@optionalarg
   \printcitestart
   \nocite{#1}%
   \@notfirstcitationfalse
   \@for \@citation :=#1\do
   {%
      \expandafter\@onecitation\@citation\@@
   }%
   \ifx\empty\@citenotetext\else
      \printcitenote{\@citenotetext}%
   \fi
   \printcitefinish
}%
\def\@onecitation#1\@@{%
   \if@notfirstcitation
      \printbetweencitations
   \fi
   \expandafter \ifx \csname\@citelabel{#1}\endcsname \relax
      \if@citewarning
         \message{\@linenumber Undefined citation `#1'.}%
      \fi
      \expandafter\gdef\csname\@citelabel{#1}\endcsname{%
\strut
\vadjust{\vskip-\dp\strutbox
\vbox to 0pt{\vss\parindent0cm \leftskip=\hsize 
\advance\leftskip3mm
\advance\hsize 4cm\strut\openup-4pt 
\rightskip 0cm plus 1cm minus 0.5cm ?  #1 ?\strut}}
         {\tt
            \escapechar = -1
            \nobreak\hskip0pt
            \expandafter\string\csname#1\endcsname
            \nobreak\hskip0pt
         }%
      }%
   \fi
   \csname\@citelabel{#1}\endcsname
   \@notfirstcitationtrue
}%
%
%
\def\@citelabel#1{b@#1}%
%
%
\def\@citedef#1#2{\expandafter\gdef\csname\@citelabel{#1}\endcsname{#2}}%
%
%
%
\def\@readbblfile{%
   \ifx\@itemnum\@undefined
      \@innernewcount\@itemnum
   \fi
   \begingroup
      \def\begin##1##2{%
         \setbox0 = \hbox{\biblabelcontents{##2}}%
         \biblabelwidth = \wd0
      }%
      \def\end##1{}
      %
      %
      \@itemnum = 0
      \def\bibitem{\@getoptionalarg\@bibitem}%
      \def\@bibitem{%
         \ifx\@optionalarg\empty
            \expandafter\@numberedbibitem
         \else
            \expandafter\@alphabibitem
         \fi
      }%
      \def\@alphabibitem##1{%
         \expandafter \xdef\csname\@citelabel{##1}\endcsname {\@optionalarg}%
         \ifx\biblabelprecontents\@undefined
            \let\biblabelprecontents = \relax
         \fi
         \ifx\biblabelpostcontents\@undefined
            \let\biblabelpostcontents = \hss
         \fi
         \@finishbibitem{##1}%
      }%
      \def\@numberedbibitem##1{%
         \advance\@itemnum by 1
         \expandafter \xdef\csname\@citelabel{##1}\endcsname{\number\@itemnum}%
         \ifx\biblabelprecontents\@undefined
            \let\biblabelprecontents = \hss
         \fi
         \ifx\biblabelpostcontents\@undefined
            \let\biblabelpostcontents = \relax
         \fi
         \@finishbibitem{##1}%
      }%
      \def\@finishbibitem##1{%
         \biblabelprint{\csname\@citelabel{##1}\endcsname}%
         \@writeaux{\string\@citedef{##1}{\csname\@citelabel{##1}\endcsname}}%
         \ignorespaces
      }%
      %
      %
      \let\em = \bblem
      \let\newblock = \bblnewblock
      \let\sc = \bblsc
      \frenchspacing
      \clubpenalty = 4000 \widowpenalty = 4000
      \tolerance = 10000 \hfuzz = .5pt
      \everypar = {\hangindent = \biblabelwidth
                      \advance\hangindent by \biblabelextraspace}%
      \bblrm
      \parskip = 1.5ex plus .5ex minus .5ex
      \biblabelextraspace = .5em
      \bblhook
      \input \bblfilebasename.bbl
   \endgroup
}%
%
%
\@innernewdimen\biblabelwidth
\@innernewdimen\biblabelextraspace
%
%
%
\def\biblabelprint#1{%
   \noindent
   \hbox to \biblabelwidth{%
      \biblabelprecontents
      \biblabelcontents{#1}%
      \biblabelpostcontents
   }%
   \kern\biblabelextraspace
}%
%
%
%
\def\biblabelcontents#1{{\bblrm [#1]}}%
%
%
\def\bblrm{\rm}%
%
%
\def\bblem{\it}%
%
%
\def\bblsc{\ifx\@scfont\@undefined
              \font\@scfont = cmcsc10
           \fi
           \@scfont
}%
%
%
\def\bblnewblock{\hskip .11em plus .33em minus .07em }%
%
%
\let\bblhook = \empty
%
%
%
\def\printcitestart{[}
\def\printcitefinish{]}
\def\printbetweencitations{, }
\def\printcitenote#1{, #1}
%
%
%
\let\citation = \@gobble
%
%
%
\@innernewcount\@numparams
%
%
\def\newcommand#1{%
   \def\@commandname{#1}%
   \@getoptionalarg\@continuenewcommand
}%
%
%
\def\@continuenewcommand{%
   \@numparams = \ifx\@optionalarg\empty 0\else\@optionalarg \fi \relax
   \@newcommand
}%
%
%
\def\@newcommand#1{%
   \def\@startdef{\expandafter\edef\@commandname}%
   \ifnum\@numparams=0
      \let\@paramdef = \empty
   \else
      \ifnum\@numparams>9
         \errmessage{\the\@numparams\space is too many parameters}%
      \else
         \ifnum\@numparams<0
            \errmessage{\the\@numparams\space is too few parameters}%
         \else
            \edef\@paramdef{%
               \ifcase\@numparams
                  \empty  No arguments.
               \or ####1%
               \or ####1####2%
               \or ####1####2####3%
               \or ####1####2####3####4%
               \or ####1####2####3####4####5%
               \or ####1####2####3####4####5####6%
               \or ####1####2####3####4####5####6####7%
               \or ####1####2####3####4####5####6####7####8%
               \or ####1####2####3####4####5####6####7####8####9%
               \fi
            }%
         \fi
      \fi
   \fi
   \expandafter\@startdef\@paramdef{#1}%
}%
%
%
%
%
\def\@readauxfile{%
   \if@auxfiledone \else 
      \global\@auxfiledonetrue
      \@testfileexistence{aux}%
      \if@fileexists
         \begingroup
            \endlinechar = -1
            \catcode`@ = 11
            \input \jobname.aux
         \endgroup
      \else
         \message{\@undefinedmessage}%
         \global\@citewarningfalse
      \fi
      \immediate\openout\@auxfile = \jobname.aux
   \fi
}%
%
%
\newif\if@auxfiledone
\ifx\noauxfile\@undefined \else \@auxfiledonetrue\fi
%
%
%
%
\@innernewwrite\@auxfile
\def\@writeaux#1{\ifx\noauxfile\@undefined \write\@auxfile{#1}\fi}%
%
%
%
\ifx\@undefinedmessage\@undefined
   \def\@undefinedmessage{No .aux file; I won't give you warnings about
                          undefined citations.}%
\fi
%
%
\@innernewif\if@citewarning
\ifx\noauxfile\@undefined \@citewarningtrue\fi
%
%
%
\catcode`@ = \@oldatcatcode
 \newpage

\head {\S0 Introduction} \endhead
\resetall
\bigskip 

Shelah Spencer \cite{ShSp:304} proved the $0-1$ law for the random graphs
$G(n,p_n)$, \newline
$p_n=n^{- \alpha}$,
$\alpha \in (0,1)$ irrational (set of nodes in 
$[n] = \{1,\dotsc,n\}$, the edges
are drawn independently, probability of edge is $p_n$).  
One may wonder what can we say on sentences
$\psi$ for which Prob$(G(n,p_n) \models \psi)$ converge to zero, Lynch 
\cite{L} asked the question and did the analysis, getting (for every $\psi$):
\medskip
\roster
\item "{$(\alpha)$}"  Prob$[G(n,p_n) \models \psi] = cn^{- \beta} + 
O(n^{- \beta-\varepsilon})$ for some $\varepsilon$ such that $\beta > 
\varepsilon > 0$
\endroster
\smallskip

\noindent
or
\roster
\item "{$(\beta)$}"  Prob$(G(n,p_n) \models \psi) = O(n^{-\varepsilon})$ for
every $\varepsilon > 0$.
\endroster
\medskip

\noindent

Lynch conjectured that in case $(\beta)$ we have
\medskip
\roster
\item "{$(\beta^+)$}"  Prob$(G(n,p_n) \models \psi) = O
(e^{-n^\varepsilon})$ for some $\varepsilon  > 0$.
\endroster
\medskip

\noindent
We prove it here.
\bigskip

\noindent
\underbar{Notation}  Let $\ell,m,n,k$ be natural numbers. \newline
Let $\varepsilon,\zeta,\alpha,\beta,\gamma$ be positive reals. \newline
$[n] = \{1,\dotsc,n\}$. \newline
$\Bbb R$ is the set of reals. \newline
$\Bbb R^+$ is the set of reals $> 0$.
\newpage

\head {\S1} \endhead
\resetall
\bigskip

\proclaim{\stag{1} Theorem}  1) For any 
first order sentence $\psi$ in the language
of graphs and irrational $\alpha \in (0,1)_{\Bbb R}$ we have (where
$p_n = n^{- \alpha}$ and Prob$(G_{n,p_n} \models \psi) \rightarrow 0$):

$$
\text{\underbar{either} } \text{Prob}(G_{n,p_n} \models \psi) \text{ is }
cn^{- \beta} + O(n^{\beta-\varepsilon}) \text{ for some reals } \beta > 
\varepsilon > 0 \text{ and } c > 0
$$

$$
\text{\underbar{or} } \text{Prob}(G_{n,p_n} \models \psi) \text{ is }
O(e^{-n^\varepsilon}) \text{ for some real } \varepsilon > 0.
$$
\medskip

\noindent
2)  However, this is not recursive.
\endproclaim
\bigskip

\demo{Proof}  We change the context generalizing it.
\enddemo
\bigskip

\definition{\stag{2} Definition of the Probability Context}
\medskip
\roster
\item "{(a)}"  $Q_n \subseteq \{ 1,\dotsc,n\},G^*_{Q_n}$ a graph on $Q_n$.
\item "{(b)}"  We consider first order sentences or formulas with vocabulary
\newline
$\subseteq \tau = \{=,Q,R\}$,
($=$ is equality, $Q$ is a monadic predicate, $R$ is a symmetric 
irreflexive binary
relation (will be ``being an edge").)
\item "{(c)}"  $G = G_{n,p_n}[G^*_Q]$ a graph on $[n]$,
$G \restriction Q = G^*_Q$,
and except this, $G$ is random with edge probability $p_n$ (i.e. for 
every edge
not included in $Q$ we flip a coin with probability $p_n$ and do it
independently for the set of edges).  We consider $G$ a $\tau$-model 
with $Q^M = Q,R$ the edge relation.
\endroster
\enddefinition
\bigskip

\remark{Remark}  The point is that $|Q|$ will be required to be just 
$< n^\varepsilon$ not say $< \text{ log}(n)$.
\endremark
\bigskip

\demo{Proof}  We consider only graphs $H$ in $\{H:H$ a graph whose set of
nodes include $Q$, moreover $H \restriction Q = G^*_Q\}$.  First, 
we repeat the proof in Shelah Spencer \cite{ShSp:304},
section 4, starting in p.105.  In our context we define
``$[H_0,H_1)$ has type $(v,e)$", it holds if $v = |H_1 \backslash H_0 
\backslash Q|$, and

$$
e = \biggl| \biggl\{ \{ x,y \} \in E(G_{n,p}):\{ x,y \} \subseteq H_1 \cup Q,
\{ x,y \} \nsubseteq H_0 \cup Q \biggr\} \biggl|,
$$

\noindent
(where for a graph $G$, $E(G)$ is the set of edges of $G$). \newline
Then define dense, sparse, safe, rigid, hinged as there adding ``over $Q$
and/or inside $G$" for definiteness.  We also define $cl_\ell(H_0;H_1)$ as in p.107, line 7.
Later we write $cl_\ell(H_0;Q)$.  All claims hold, but
arriving to Theorem \scite{3} (bottom of p.107) we should be careful.
We consider only embeddings which are the identity on $Q$.
\enddemo
\bigskip

\proclaim{\stag{3} Lemma}  1) Let $\ell^* \in \Bbb N$.
For every small enough $\varepsilon > 0$, for some $\xi > 0$, for every
$n$ large enough, if $|Q| \le n^\xi,Q \subseteq [n]$ we have:
if $(H_0,H_1)$ is safe of type $(v,e)$ and $f$ embeds $H_0$ into $G$ (and
$f$ is the identity on $Q$) and $|H_1 \backslash Q| \le \ell^*$ 
\underbar{then}:

$$
\text{Prob} \biggl( \neg[n^{v-\alpha \, e-\varepsilon} < N(f,H_0,H_1) <
n^{v-\alpha \, e+\varepsilon} ] \biggr) < e^{-n^\xi}
$$

\noindent
(where $N(f,H_0,H_1)$ is the number of extensions $g:H_1 \rightarrow G$
satisfying: \newline
$x \in H_0 \Rightarrow g(x) = f(n)$ and $\{x,b\} \in E(H_1),
b \notin H_0 \Rightarrow \{g(x),g(y)\} \in E(G)$). \newline
2)  Let $\varepsilon \in \Bbb R^+$ and $\ell^* \in \Bbb N$ be given, 
then for some $\xi > 0$ for every $n$ large enough and any
$Q \subseteq [n],|Q| \le n^\varepsilon$ and graph $G^*_Q$ on $Q$ we 
consider only embeddings which are the identity on $Q$.  Then
\medskip
\roster
\item "{$(*)$}"  if $H_1$ is a graph with $|H_1 \backslash Q| \le \ell^*$,
$H_0 \subseteq H_1$, we assume $f$ embeds $H_0$ into $Q$,
$f$ is the identity on $H_0$ and $(H_0,H_1)$ is rigid
then:
$$
\text{Prob} \biggl( N(f,H_0,H_1,G_{n,p_n}) > 0 \biggr) < n^{- \varepsilon}.
$$
\endroster
\endproclaim
\bigskip

\demo{Proof}  1) As in \cite[Theorem 3,p.107]{ShSp:304} + extra computation 
by the central limit theorem \underbar{or} see \cite[\S5]{Sh:550} for 
more. \newline
2)  As in \cite{ShSp:304}.
\enddemo
\bigskip

\proclaim{4 Lemma}  For any $k,m \in \Bbb N$ there are $\ell^*$ and
$\varepsilon^* > 0$ depending on $k$ only such that the following holds:
\medskip
\roster
%
\item "{$(*)$}"  For any formula $\psi = \psi(x_1,\dotsc,x_m)$ of quantifier
depth $\le k$ in the vocabulary $\{=,Q,R\}$ there is a formula $\theta_\psi 
= \theta_\psi(x_1,\dotsc,x_m)$ in the vocabulary \newline
$\{=,Q,R\}$ such that:
\item "{$(**)$}"  For every $n$ large enough, $Q \subseteq \{1,\dotsc,n\},
|Q| \le n^{\varepsilon^*}$, and graph $G^*_Q$ on $Q$ and
$G = G_{n,p_n}[G^*_Q]$ such that the small probability cases from Lemma
\scite{3(1)},\scite{3(2)} (for 
$(H_1,H_2)$ of type $(v,e),v \le 2 \ell^*$), or just
$\otimes^1_{\ell^*} + \otimes^2_{\ell^*}$ below do not occur, we have:
\item "{$(***)$}"  for every $a_1,\dotsc,a_m \in \{1,\dotsc,m\}$ we have
$$
\gather
(\{1,\dotsc,n\},Q,R) \models \psi[a_1,\dotsc,a_m] \, \text{ \underbar{iff}}\\
\biggl( Q \cup \{a_1,\dotsc,a_m\},Q,R \restriction (Q \cup \{ a_1,
\dotsc,a_m\}) \biggr) \models \theta_\psi[a_1,\dotsc,a_m]. \endgather
$$
\endroster
\medskip

\noindent
where
\roster
\item "{$\bigotimes^1_{\ell^*}$}"  if $(H_0,H_1)$ is safe (so $Q \subseteq
H_0$) \newline
$|H_1 \backslash Q| \le \ell^*,H_0 \subseteq G_{n,p_n}[G^*_Q]$ then we can
extend $\text{id}_{H_0}$ to an embedding $g$ of $H_1$ into 
$G_{n,p_n}[G^*_Q]$ such that \newline
$c \ell_{\ell^*} \left( g(H_1),G_{n,p}[G^*_q] \right)
= g(H_1) \cup c \ell_{\ell^*} \left( f(H_0,G_{n,p_n}[G^*_Q] \right)$
\medskip
\noindent
\item "{$\bigotimes^2_{\ell^*}$}"  if $(H_0,H_1)$ is rigid, $|H_1 \backslash
Q| \le \ell^*,H_0 = G^*_Q$ then there is no extension of $f$ of
$\text{id}_{H_0}$ to an embedding of $H_1$ into $G_{n,p_n}[G^*_Q]$.
\endroster
\endproclaim
\bigskip

\demo{Proof}  Similar to the proof in \cite{ShSp:304}, 
(and is a particular case of \cite[\S2]{Sh:467} (see related)).
\enddemo
\bigskip

\demo{Proof of Theorem \stag{1}}  Part (1)  Let $\theta_\psi$ 
be from the analysis
(i.e. Lemma \scite{4} for the $\psi$ from Theorem \scite{1}) 
for the original sentence $\psi$.
\enddemo
\bigskip

\subhead {Case A} \endsubhead  For some finite graph $G^*$ on say 
$\{1,\dotsc,m^*\}$ we have $G^* \models \theta_\psi$.

In this case the probability that $G^*$ can be embedded into $G_{n,p_n}$ is
$\ge O(n^{- \beta})$ for some $\beta \in (0,\infty)$ if $n \ge m^*$ of 
course; so this
means that one of the $\le n^{m^*}$ possible mapping is an embedding, but more
convenient is to consider the event $G \restriction [m^*] = G^*$ which also
has probability $\ge n^{- \beta}$ for some $\beta$.  Now modulo this event
the probability that the conclusion of Lemma \scite{4} fails is (for $n$ large
enough) much smaller than $n^{-m^*}$.  So we can assume that for 
$G \restriction [m^*] \cong G^*$ and that the conclusion of Lemma \scite{4}
holds for this.  Now check and if we succeed by Lemma \scite{4}, 
we are done, i.e. the
probability that $G_{n,p_n} \models \psi$ is quite high.
\bigskip

\subhead{Case B} \endsubhead  For no finite graph $G^*,G^* \models
\theta_\psi$. \newline

Choose $\ell^* \in \Bbb N$ large enough as needed for our sentence $\psi$
in Lemma 4.

Let $\zeta \in \Bbb R^+$ be such that: \newline
$v \in \{ 0,\dotsc,2 \ell^*\}, e \in \Bbb N \Rightarrow |v - \alpha e| \ge
\zeta$ and it satisfies the requirements on $\zeta$ in Lemma \scite{3(2)} (for
$2 \ell^*$ (readily follows).) \newline
(The $2 \ell^*$ rather than $\ell^*$ is for the bound on Prob$({\Cal E}_2)$.)
Clearly $\zeta$ exists and if $(H_0,H_1)$ is hinged and $|H_1 \backslash
H_0| \le \ell^*$ and $(H_0,H_1)$ is of type $(v,e)$ then $v - \alpha e <
- \zeta$.
\medskip

\noindent
Let $\varepsilon(\ell^*),\xi$ be such that:
\medskip
\roster
\item "{$(a)$}"  $\varepsilon(\ell^*) \in \Bbb R^+$ and
$\varepsilon(\ell^*) < \zeta/2,\xi < \zeta/2$
\item "{$(b)$}"  in Lemma \scite{3(1)} $\varepsilon(\ell^*),\xi$ satisfies the
requirements of $\varepsilon,\xi$ respectively.
\endroster
\medskip

\noindent
We shall prove that for $n$ large enough Prob$(G_{n,p_n} \models \psi)$
is $\le e^{-(n^\xi)}$, this is enough.

For any $G = G_{n,p_n}$, we define by induction on $j \le n$, a subset
$P_j = P_j[G]$ of $\{1,\dotsc,n\}$ as follows:

$$
P_0 = \emptyset
$$

$$
\align
P_{j+1} = P_j \cup \{H:&P_j \subseteq H \subseteq G,|H \backslash P_j|
\le \ell^*,H \ne P_j \text{ and} \\
 &(P_j,H) \text{ is rigid in } G \}.
\endalign
$$

\noindent
For some $j(*) < n$ we have $P_{j(*)} = P_{j(*)+1}$ (hence
$P_{j(*)+1} = P_{j(*)+2}$, etc). \newline
If $|P_{j(*)}| \le n^{\varepsilon(\ell^*)}$ and $\otimes^1_{\ell^*}$ holds
then, (as $P_{j(*)} = P_{j(*)+1}$) this implies $\otimes^2_{\ell^*}$ and
then by Lemma \scite{4} we are done ($P_{j(*)}$ is $Q$).  So
it is enough to give an upper bound of the form $e^{-n^\varepsilon}$ to the
probability Prob$({\Cal E}_1) + \text{ Prob}({\Cal E}_2)$ were 
${\Cal E}_1$ is the event $|P_{j(*)}| > n^{\varepsilon(\ell^*)}$ and 
${\Cal E}_2$ is the event $|P_{j(*)}| \le n^{\varepsilon(\ell^*)} 
\and [\otimes^1_{\ell^*}$ fails]. \newline
\medskip
\noindent
\underbar{On Prob$({\Cal E}_1)$}.  If $|P_{j(*)}| \ge n^{\varepsilon(\ell^*)}$
then we can find
$a_{j,\ell}$ for $j < [n^{\varepsilon(\ell^*)} / \ell^*]$ and 
$\ell < \ell_j \le \ell^*$ such that \newline
$\biggl( H_i \cap \{ a_{i,\ell}:\ell < \ell_i\},\{a_{i,\ell}:\ell < 
\ell_i\} \biggr)$ (in $G$) is rigid of type $(v_i,e_i)$ where \newline
$H_i =: \{ a_{j,\ell}:j < i \text{ and } \ell < \ell_j \}$ 
(so we may have not used all
$P_{j(*)}$).  Clearly there is a real
$\zeta > 0$ depending on $\ell^*,\alpha$ only such that 
$v_i - e_i \alpha \le -\zeta$,
(simply, there are only finitely many possible pairs $(v,e)$).
\medskip

Let $I$ be a sequence describing
this situation, i.e. it contains

$$
\langle \ell_i:i < [n^{\varepsilon(\ell^*)}/\ell^*] \rangle
$$

$$
\{((i_1,m_1),(i_2,m_2)):a_{\ell_1,m_1} = a_{i_2,m_2}\}
$$ 

$$
\{(i,m_1,m_2):a_{i,m_1}R^Ga_{i,m_2}\}.
$$
\medskip

\noindent
There are 
$\dsize \prod_{i < [n^{\varepsilon(\ell^*)}/\ell^*]}
(\ell^* \times (\ell^* \times
i)^{\ell^*} \times 2^{2 \ell^*})$ possible such sequences $I$ (an overkill).
\newline

\noindent
[Why?  The ith term in the product is an upper bound on the number of choices
in stage $i$, there $\ell^*$ is the number of possible $\ell_i,\ell^* 
\times i$ is
an upper bound on the number $|\{a_{j,\ell}:j < i,\ell < \ell_j\}|$,
$(\ell^* \times i)^{\ell^*}$ is an upper bound to the number of choices of
$\langle a_{i,\ell}:\ell < \ell^*,a_{i,\ell} \in \{ a_{j,s}:j < i,
s < \ell_j\} \rangle$, and $2^{2 \ell^*}$ is an upper bound to the number of
possible $G \restriction \{ a_{i,\ell}:\ell < \ell_i\}$].

Now for some constants $c_0,c_1$ depending only on $\ell^*$ (i.e. $\psi$)
this number is $\le c^{n^{\varepsilon(\ell^*)}/\ell^*}_0 \times
[(n^{\varepsilon(k^*)}/\ell^*)!]^{\ell^*} \le
n^{\varepsilon(\ell^*)n^{\varepsilon(\ell^*)}}$.
For each $I$ the number of
possibilities for the $a_{i,\ell}$ is $\le \dsize \prod_i n^{v_i}$, and
the probability it holds in $G$ is $\dsize \prod_i n^{- \alpha e_i}$, 
hence the expected value is

$$
\dsize \prod_i n^{(v_i - \alpha e_i)} \le \dsize \prod_i n^{-\zeta} =
n^{-\zeta(n^{\varepsilon(\ell^*)}/\ell^*)}.
$$

\noindent
So the expected number of number of such $\langle a_{i,\ell}:i <
n^{\varepsilon(\ell^*)}/\ell^* \text{ and } \ell < \ell_i \rangle$ \newline
for some $I$ is $\le n^{(\varepsilon(\ell^*)-\zeta)n^{\varepsilon(\ell^*)}}$ 
and 
as we have $\varepsilon(\ell^*) < \zeta$ the conclusion should be clear.
\bigskip

\noindent
\underbar{Probability of ${\Cal E}_2$}.  Should be 
clear by Lemma \scite{3(1)}; i.e.
except suitably small probability the number of extensions of $f$ to
embedding of $H_1$ is much larger than the number of such extensions failing
the requirement in $\otimes^1_{\ell^*}$.
\bigskip

\demo{Proof of Theorem \stag{1}-part (2)}  In non-trivial 
cases for some $\ell$ and pair \newline
$(H_0,H_1)$ we have $H_1 \ne H_0$ and $H_1 \subseteq cl_\ell(H_0)$. \newline
  Now for $n$ large enough
(if $|cl_\ell(H_0)| \ll \text{ log } n$), \newline
on $cl_\ell(H_0)$ in $G_{n,p_n}$ we can
interpret arithmetic on $cl_\ell(H_0)$ (with parameters) and all subsets and
all second place relations.  Fix $H_0,\ell$. \newline

For a sentence $\psi$ speaking on $\Bbb N \restriction k$,
(or $2^k$) we can compute
$\psi^*$ in the vocabulary of graphs saying
\medskip
\roster
\item "{$(*)$}"   there is a copy $H'_0$ of $H_0$ such that
$$
\Bbb N \restriction |cl_\ell(H'_0) \models \psi^*.
$$
\endroster
\medskip

\noindent
So for every function $h:\Bbb N \rightarrow \Bbb N$ converging to infinity

$$
\text{Lim inf}_n 
\biggl(\text{Prob}(G_{n,p_n} \models \psi^*)/n^{-h(n)} \biggr)
\ge 1 \text{\underbar{ iff} }
\dsize \bigvee_k[\Bbb N \restriction k \models \psi].
$$

\noindent
But the set $\{ \psi:(\exists k)[|N \restriction k| \models \psi] \}$ is 
like the
set of sentences having a finite model (i.e. same Turing degree)
so is not recursive.
\enddemo
\bigskip

\remark{Concluding Remarks}  1) In fact, we have to consider $P_j$
(in case $B$ during the proof of Theorem 1) only for $j \le 2^r$, where $r$
is the quantifier depth of the sentence $\psi$ (for which we are proving
Theorem \scite{1}). 
\endremark
\newpage

REFERENCES.  
\bibliographystyle{lit-plain}
\bibliography{lista,listb,listx}

\enddocument

\bye